\documentclass{article}

\usepackage[T1]{fontenc}
\usepackage[latin1]{inputenc}

\usepackage{amsmath,amsthm,amssymb,mathtools,array,centernot,paralist}

\newtheorem{theorem}{Theorem}[section]
\newtheorem{lemma}[theorem]{Lemma}
\newtheorem{corollary}[theorem]{Corollary}
\theoremstyle{definition}
\newtheorem{example}{Example}[section]
\newtheorem*{remark}{Remark}

\newcommand{\Z}{\ensuremath{\mathbb{Z}}}

\DeclareMathAlphabet{\matris}{T1}{cmss}{m}{sl}
\newcommand{\mat}[1]{\ensuremath{\matris{#1}}}
\newcommand{\vek}[1]{\ensuremath{\boldsymbol{#1}}}

\DeclareMathOperator{\adj}{adj}

\newenvironment{ekvsystem}[1][3]{%
  \left\{
    \begin{array}{*{#1}{@{}r@{}>{\;}c<{\;}}@{}r}
    }{%
    \end{array}
  \right.
}

\usepackage{color,soul,framed}
\sethlcolor{yellow}
\definecolor{shadecolor}{cmyk}{0,0,1,0}

\title{Number of Solutions of\\ Linear Congruence Systems}
\author{%
  Marcus Nilsson\\
  \texttt{marcus.nilsson@lnu.se}\\
  Linnaeus University\\
  Sweden
  \and
  Robert Nyqvist\\
  \texttt{robert.nyqvist@bth.se}\\
  Blekinge Institute of Technology\\
  Sweden
}

\begin{document}

\maketitle

\begin{abstract}
\end{abstract}

\section{Introduction}	
We will consider the system of linear congruenses, 
\begin{equation}
  \label{eq:lcs}
  \left\{
    \begin{array}{@{}l@{\;}c@{\;}l}
      a_{11} x_1 + a_{12} x_2 + \cdots + a_{1n} x_n & \equiv & b_1 \pmod{m} \\
      a_{21} x_1 + a_{22} x_2 + \cdots + a_{2n} x_n & \equiv & b_2 \pmod{m} \\
      & \vdots & \\
      a_{n1} x_1 + a_{n2} x_2 + \cdots + a_{nn} x_n & \equiv & b_n \pmod{m},
    \end{array}
  \right.
\end{equation}
where $m$, $n$ are positive integers, and all $a_{ij}$, $b_i$ are
integers. We are interested in finding an expression for the number of
solutions to this system. Let $\mat{A}=(a_{ij})$ denote the coefficient
matrix. It is well known that this system has a unique solution if and only if
$\det(\mat{A})$ and $m$ are relatively prime, see for
example~\cite{Rosen2011}. But, how many solutions do we have if
$\det(\mat{A})$ and $m$ have a common divisor greater than 1? The authors
became interested in this question when they where working on
multidimensional $p$-adic monomial dynamical systems,
see~\cite{NilsNyqv2010}. The problem was solved by Butson and
Stewart~\cite{ButSte1955}, by rewriting the system into Smith normal form.
The solvability of linear congruence systems and algorithms for finding a
solutions have been of interest for many mathematicians and computer
scientists over the years.  Example of other contributions to the problem on
solving systems of linear congruences can be found
in~\cite{Cul1901,DolzStur2001-1,DolzStur2001-2,McCa1975,Smar1987}.

Compared to Butson and Stewart we will use a more direct method. We use
Gaussian elimination with successive reduction and the Chinese remainder
theorem instead of the Smith normal form. We will find a different formula
for the number of solutions than in~\cite{ButSte1955}. The algorithm used in
the proof of the formula can with small changes also be used to find all
incongruent solutions to the system, since we only have used elementary
methods.

The paper is organized as follows: In Section \ref{sec:notations} we give
definitions and notations together with some theorems about solvability of
systems of linear congruences.  Section \ref{sec:homogeneous} is the main
section in which we derive the formula for the number of solutions to
homogeneous systems. Applications to inhomogeneous systems are mentioned in
Section \ref{sec:inhomogeneous}. In Appendix~\ref{sec:algorithm} we present an
algorithm in pseudo code for calculating the number of solutions. The
algorithm also gives us all the solutions. In Section \ref{sec:discussion} we
discuss generalizations of our results.

\section{Notations}
\label{sec:notations}
Two solutions of the linear congruence system~\eqref{eq:lcs} are said to be
\emph{incongruent modulo $m$} if they differ at least in one coordinate
modulo $m$.  We want to find the number of incongruent solutions modulo $m$.
Let \mat{A} denote the coefficient matrix, \vek{x} the
vector of the indeterminates, and \vek{b} the vector of the elements on the
right hand side in the system, that is,~$\mat{A} = (a_{ij})$, $\vek{x} =
(x_1, x_2, \ldots, x_n)$ and $\vek{b} = (b_1, b_2, \ldots, b_n)$.  Then the
system~\eqref{eq:lcs} can be written in matrix form as
\begin{equation}
  \label{eq:matrixform}
  \mat{A}\vek{x} \equiv \vek{b} \pmod{m}.
\end{equation}
Let $\eta(\mat{A}, \vek{b}, m)$ denote the number of incongruent
solutions modulo $m$ to the congruence~\eqref{eq:matrixform}.

Let $\adj(\mat{A})$ denote the adjoint matrix of a square matrix \mat{A} of
order~$n$.  It is known from linear algebra that $\adj(\mat{A})$ has the
following properties:
\[
  \det(\adj(\mat{A})) = \det(\mat{A})^{n-1}
\]
and
\[
  \mat{A}\, \adj(\mat{A})
  =
  \adj(\mat{A})\, \mat{A}
  =
  \det(\mat{A})\, \mat{I},
\]
where $\mat{I}$ is the identity matrix.

\begin{theorem}
  \label{th:estimate}
  Let \mat{A}, \vek{b} and $m$ be as above.
  Then
\[
  \eta(\mat{A}, \vek{b}, m) \leq (\det(\mat{A}), m)^n.
\]
If $(\det(\mat{A}), m) = 1$, then $\eta(\mat{A}, \vek{b}, m) = 1$.
\end{theorem}
\begin{proof}
  We multiply the congruence system $\mat{A} \vek{x} \equiv
  \vek{b} \pmod{m}$ by $\adj(\mat{A})$ and get
  \begin{gather}
    \adj(\mat{A}) \, \mat{A} \vek{x}
    \equiv \adj(\mat{A}) \, \vek{b} \pmod{m} \notag\\
    \Leftrightarrow \notag\\
    \det(\mat{A}) \, \vek{x}
    \equiv \adj(\mat{A}) \, \vek{b} \pmod{m},
    \label{eq:adj}
  \end{gather}
  which is solvable if and only if $\det(\mat{A})$ divides all the elements in
  $\adj(\mat{A}) \vek{b}$.  If that is the case then this system
  has $(\det(\mat{A}), m)^n$ different solutions.  If that is not the case,
  then the original system has no solutions, since any solutions is also
  solutions to the rewritten system~\eqref{eq:adj}.  When we multiply by
  the matrix~$\adj(\mat{A})$ it might happen that we introduce new solutions.
  This proves that the inequality $\eta(\mat{A}, \vek{b}, m) \leq
  (\det(\mat{A}), m)^n$ holds.

  Assume that $(\det(\mat{A}), m) = 1$.  From~\eqref{eq:adj} it follows
  that
  \[
    \vek{x}
    \equiv
    \det(\mat{A})^{-1} \, \adj(\mat{A}) \, \vek{b}
    \pmod{m}
  \]
  is a solution to the original system. Hence $\eta(\mat{A},\vek{b},m)=1$.
\end{proof}

\begin{theorem}
  \label{thm:1}
  Let $\mat{A}$, $\vek{b}$ and $m$ be as above.  Assume that $m = m_1
  \cdots m_k$, where the integers $m_1, \ldots, m_k$ are pairwise relatively
  prime.  Then
  \[
    \eta(\mat{A}, \vek{b}, m)
    =
    \eta(\mat{A}, \vek{b}, m_1)
    \cdots
    \eta(\mat{A}, \vek{b}, m_k).
  \]
  Hence, $\eta$ is multiplicative with respect to $m$.
\end{theorem}
\begin{proof}
  Since $m_i \mid m$ for each $i$, any solution to
  $\mat{A} \vek{x} \equiv \vek{b} \pmod{m}$ is also a
  solution to every $\mat{A} \vek{x} \equiv \vek{b} \pmod{m_i}$.
  Assume that $\vek{x}_i = (x_{i1}, x_{i2}, \ldots,
  x_{in})$ is a solution to the congruence system $\mat{A} \vek{x} \equiv
  \vek{b} \pmod{m_i}$, for $i = 1, 2, \ldots, k$.  Then for any $k$-tuple
  $(\vek{x}_1, \vek{x}_2, \ldots, \vek{x}_k)$ of solutions of the systems
  modulo $m_i$ we construct a solution $\vek{x} = (x_1, x_2, \ldots, x_n)$
  modulo $m$, by defining $x_i$ to be the unique solution, due to the Chinese
  Remainder Theorem, modulo~$m$, to the system
  \[
    \begin{ekvsystem}[1]
      x & \equiv & x_{1i} \pmod{m_1} \\
      x & \equiv & x_{2i} \pmod{m_2} \\
      & \vdots & \\
      x & \equiv & x_{ki} \pmod{m_k}\rlap{.}
    \end{ekvsystem},
  \]
  Since there are $\eta(\mat{A}, \vek{b}, m_1) \cdots \eta(\mat{A}, \vek{b},
  m_k)$ possible such $k$-tuples of solutions modulo~$m_i$, the theorem
  follows (two different $k$-tuples can not generate the same solution modulo
  $m$).
\end{proof}

\section{Homogeneous Systems}
\label{sec:homogeneous}

According to Theorem~\ref{thm:1} we only have to consider congruence
 systems
modulo a prime power.  So, from now on $m = p^k$,
where $p$ is a prime number and $k$ a positive integer.  Let $d = (\mat{A},
p^k)$ be the greatest common divisor of all the
elements in the matrix $\mat{A}$ and the integer $p^k$, that is,
\[
  (\mat{A}, p^k) = (a_{11}, a_{12}, \ldots, a_{nn}, p^k).
\]
Hence, $d=p^e$ for some non-negative integer $e \leq k$.

\begin{lemma}
  \label{lma:homogen1}
  Let $d=(\mat{A}, p^k) = p^e$ and let $\mat{R} = \mat{A}/d$.  Then
  \[
    \eta(\mat{A}, \vek{0}, p^k)
    =
    d^n \eta(\mat{R}, \vek{0}, p^k/d)
    =
    p^{en} (\mat{R}, \vek{0}, p^{k-e}).
  \]
  Note that $\eta(\mat{R}, \vek{0}, 1) = 1$.
\end{lemma}
\begin{proof}
  Let $q = p^k/d = p^{k-e}$.  Assume that $\vek{y} \in \Z_q^n$ is a solution
  to the congruence system $\mat{R} \vek{x} \equiv \vek{0} \pmod{p^k/d}$.
  Then it follows from the definition of \mat{R} that \vek{y} is also a
  solution to the original system $\mat{A} \vek{x} \equiv \vek{0}
  \pmod{p^k}$.  This is also true for all elements on the form
  \begin{equation}
    \vek{x}
    =
    \vek{y} + q \vek{v}
    =
    \vek{y} + p^{k-e} \vek{v},\qquad
    \text{where $\vek{v} \in \Z_d^n$},
    \label{eq:all-solutions}
  \end{equation}
  since
  \[
    \mat{A}(\vek{y} + q \vek{v})
    \equiv
    \frac{p^k}{d} \mat{A} \vek{v}
    \equiv
    p^k \mat{R} \vek{v}
    \equiv
    \vek{0} \pmod{p^k}
  \]
  and all the elements in $\mat{R} = \mat{A}/d$ are integers.  The components of
  the vector~$q\vek{v}$ are all non-negative integers less than $q d = p^k$.
  Hence, two vectors of the form~\eqref{eq:all-solutions} with different
  \vek{v} is incongruent modulo~$p^k$.
    
  Next, we prove that different solutions
  modulo $q$ are lifted to different solutions modulo $p^k$.
    If $\vek{x}_1 = \vek{y}_1 + q\vek{v}_1$ and $\vek{x}_2 =
    \vek{y}_2 + q\vek{v}_2$ are two solutions to $\mat{A}\vek{x} \equiv
    \vek{0} \pmod{p^k}$ where $\vek{y}_1$ and $\vek{y}_2$ are solutions to
    $\mat{R}\vek{y} \equiv \vek{0} \pmod{q}$, then
  \[
    \text{%
        $
          \vek{x}_1 - \vek{x}_2
          =
          \vek{y}_1 - \vek{y}_2 + q (\vek{v}_1 - \vek{v}_2)
          \equiv
          \vek{y}_1  - \vek{y}_2 \pmod{q}.
        $
      }%
  \]
    Hence, if $\vek{x}_1 = \vek{x}_2$, then $\vek{y}_1$ and $\vek{y}_2$ is
    congruent modulo~$q$.
  %
  It therefore follows that
  \[
    \eta(\mat{A}, \vek{0}, p^k) \geq d^n \eta(\mat{R}, \vek{0}, p^k/d).
  \]
  It remains to prove that all the solutions of the
  original congruence system $\mat{A}\vek{x} \equiv \vek{0} \pmod{p^k}$ can
  be written on the form~\eqref{eq:all-solutions} where the corresponding
  vector~\vek{y} is a solution to $\mat{R} \vek{x} \equiv \vek{0}
  \pmod{p^k/d}$.  Assume that $\vek{x} \in \Z_{p^k}^n$ is a solution to
  $\mat{A}\vek{x} \equiv \vek{0} \pmod{p^k}$.  By the division algorithm,
  applied on each component of \vek{x}, there are vectors $\vek{y} \in
  \Z_q^n$ and $\vek{v} \in \Z_d^n$ such that $\vek{x} = \vek{y} + q \vek{v}$.
  Then
  \[
    \vek{0}
    \equiv
    \mat{A}\vek{x}
    \equiv
    \mat{A}\vek{y} + \frac{p^k}{d} \mat{A} \vek{v}
    \equiv
    d\mat{R}\vek{y} + p^k \mat{R}\vek{v}
    \equiv
    p^e \mat{R}\vek{y} \pmod{p^k},
  \]
  or equivalent
  \[
    \mat{R}\vek{y} \equiv \vek{0} \pmod{p^{k-e}}.
  \]
  This proves the theorem.
\end{proof}

Let $\mat{A}_n = \mat{A}$ and $p^{l_n} = p^k/d_n$, where $l_{n+1} = k$ and
$d_n = p^{e_n} = (\mat{A}_n, p^k)$.  Set $\mat{R}_n = \mat{A}_n/d_n$.  If $e_n = k$,
then $\mat{A}_n \equiv \mat{0} \pmod{p^k}$ and
\[
  \eta(\mat{A}_n, \vek{0}, p^k) = p^{kn},
\]
since any vector in $\Z_{p^k}^n$ is a solution to $\mat{A}_n\vek{x} \equiv
\vek{0} \pmod{p^k}$.

Assume that $e_n < k$.  Then $(\mat{R}_n, p) = 1$ and therefore we can find an
element in $\mat{R}_n$ that is relatively prime to~$p$, say~$r_{11}$ after a
possible rearrangement of rows and columns.  Note that $r_{11}$ is then
invertible modulo $p^{l_n}$.  By Gaussian elimination we can rewrite
\begin{equation}
  \mat{R}_n\vek{x} \equiv \vek{0} \pmod{p^{l_n}}
  \label{eq:Rxmodpl}
\end{equation}
to get the equivalent system
\begin{equation}
  \begin{ekvsystem}[5]
    r_{11}x_1 &+& r_{12}x_2 &+& r_{13}x_3 &+& \cdots &+& r_{1n}x_n &\equiv&
    0 \pmod{p^{l_n}}\\
    && a_{22}'x_2 &+& a_{23}'x_3 &+& \cdots &+& a_{2n}'x_n &\equiv&
    0 \pmod{p^{l_n}}\\
    && && && && & \vdots \\
    && a_{n2}'x_2 &+& a_{n3}'x_3 &+& \cdots &+& a_{nn}'x_n &\equiv&
    0 \pmod{p^{l_n}}\rlap{,}
    \end{ekvsystem}
  \label{eq:system-pl-reduced}
\end{equation}
where
\[
  a_{ij}'
  \equiv
  r_{ij} - r_{11}^{-1} r_{i1} r_{1j}
  \pmod{p^{l_n}}
\]
for $i, j = 2, 3, \ldots, n$.
Note that we do not change the number of
solutions of the system~\eqref{eq:Rxmodpl} since the two systems are equivalent
because Gaussian transform is invertible.
%
%
  Hence, the
two systems \eqref{eq:Rxmodpl} and \eqref{eq:system-pl-reduced} have equally
many solutions modulo~$p^{l_n}$ (any solution to one of the systems is also a
solution to the other).  Let $\mat{A}_{n-1} = (a_{ij}')_{2\leq i,j \leq n}$
and $\vek{x}' = (x_2, x_3, \ldots, x_n)$, see~\eqref{eq:system-pl-reduced}.
Note that the matrix~$\mat{A}_{n-1}$ depends on $\mat{A}_n$ and the choice of
$r_{11}$ in $\mat{R}_n$.  Since $r_{11}$ is invertible modulo~$p^{l_n}$, any
solution to the congruence system
\begin{equation}
  \label{eq:An-1}
  \mat{A}_{n-1} \vek{x}' \equiv \vek{0} \pmod{p^{l_n}}
\end{equation}
can be extended with respect to $x_1$ in an unique way to a solution
to~\eqref{eq:system-pl-reduced}.  Hence, the number of solution
to~\eqref{eq:An-1} is equal to the number of solutions to~\eqref{eq:Rxmodpl}.
This proves that
\begin{equation}
  \label{eq:eta-induction}
  \eta(\mat{A}_n, \vek{0}, p^k)
  =
  d_n^n \eta(\mat{R}_n, \vek{0}, p^{l_n})
  =
  d_n^n \eta(\mat{A}_{n-1}, \vek{0}, p^{l_n}),
\end{equation}
according to Lemma~\ref{lma:homogen1}.  Further, we have that
\begin{equation}
  \label{eq:det-AnR}
  \det(\mat{A}_n) = d_n^n\det(\mat{R}_n),
\end{equation}
since $\mat{A}_n = d_n \mat{R}_n$.  From~\eqref{eq:Rxmodpl}
and~\eqref{eq:system-pl-reduced} it follows that
\[
  \det(\mat{R}_n)
  \equiv
  (-1)^{1+1} r_{11} \det(\mat{A}_{n-1})
  \pmod{p^{l_n}},
\]
and therefore is
\begin{equation}
  \label{eq:det-induction}
  (\det(\mat{R}_n), p^{l_n}) = (\det(\mat{A}_{n-1}), p^{l_n}),
\end{equation}
since $(r_{11}, p) = 1$.

\begin{remark}
  Note that if we choose an integer to represent
  $r_{11}^{-1}$ which is a multiplicative inverse modulo $p^j$ to $r_{11}$,
  where $j \geq l_n$, then~\eqref{eq:det-induction} can be reformulated
  as
\begin{equation}
  \label{eq:det-induction-2}
  (\det(\mat{R}_n), p^i) = (\det(\mat{A}_{n-1}), p^i),
\end{equation}
for all integers $i$ such that $1 \leq i \leq j$.
\end{remark}

\begin{theorem}
  \label{thm:homogen-system}
  Let \mat{A} be a square matrix of order $n$ with integer elements, $p$ a
  prime number and $k$ a positive integer.  Set $\mat{A}_n = \mat{A}$.  Let
  $\mat{A}_i$ denote the matrix consisting of the $i$ last rows and
  columns of \mat{A} after the $(n-i)$th step of Gaussian elimination
  including a possible reordering of rows or columns, where $i = n - 1,
  \ldots, 3, 2$.  Further, for $i = n, \ldots, 3, 2$ define $d_i$, $e_i$ and
  $l_i$ recursively by
  \[
    p^{l_i} = \frac{p^{l_{i+1}}}{d_i}
    \quad\text{and}\quad
    d_i = p^{e_i} = (\mat{A}_i, p^{l_{i+1}})
  \]
  with $l_{n+1} = k$.  Hence, $l_i = l_{i+1} - e_i$.  Then
  \[
    \eta(\mat{A}, \vek{0}, p^k)
    =
    (\det(\mat{A}), d_2d_3^2 \cdots d_n^{n-1} p^k).
  \]
\end{theorem}
\begin{proof}
  We prove this by induction over $n$.
  First assume that $n = 2$ and set $\mat{R}_2 = \mat{A}_2/d_2$, where $d_2 =
  (\mat{A}_2, p^{l_3}) < p^{l_3}$.  Note that $k = l_3$ in this case.  The
  congruence system $\mat{R}_2\vek{x} \equiv \vek{0} \pmod{p^{l_2}}$ is given
  by
  \begin{equation}
    \label{eq:system2}
    \begin{ekvsystem}[2]
      r_{11} x_1 & + & r_{12} x_2 & \equiv & 0 \pmod{p^{l_2}} \\
      r_{21} x_1 & + & r_{22} x_2 & \equiv & 0 \pmod{p^{l_2}}\rlap{,}
    \end{ekvsystem}
  \end{equation}
  where $p^{l_2} = p^{l_3}/d_2$ and
  \[
    \mat{R}_2
    =
    \begin{pmatrix}
      r_{11} & r_{12} \\
      r_{21} & r_{22}
    \end{pmatrix}
  \]
  is an integer matrix.  Since $(\mat{R}_2, p) = (r_{11}, r_{12}, r_{21},
  r_{22}, p) = 1$ one of the matrix entries must be relatively prime to $p$.
  Assume that it is $r_{11}$.  Let $r_{11}^{-1}$ denote the multiplicative
  inverse of $r_{11}$ modulo $p^{l_2}$.  From the first congruence in the
  system~\eqref{eq:system2} we have
  \[
    x_1 \equiv  - r_{11}^{-1} r_{12} x_2 \pmod{p^{l_2}}.
  \]
  Hence, $x_1$ is uniquely determined by $x_2$ modulo $p^{l_2}$.  By putting
  this expression for $x_1$ in the second congruence in the
  system~\eqref{eq:system2} we get
  \[
    (r_{22} - r_{11}^{-1} r_{12} r_{21}) x_2 \equiv 0 \pmod{p^{l_2}},
  \]
  or equivalent
  \[
    (r_{11} r_{22} - r_{12} r_{21}) x_2 \equiv 0 \pmod{p^{l_2}}.
  \]
  The coefficient for $x_2$ is $\det(\mat{R}_2)$ and from the theory of linear
  congruences we know that this equation has $(\det(\mat{R}_2),p^{l_2})$
  incongruent solutions modulo $p^{l_2}$.  We conclude that $\eta(\mat{R}_2,
  \vek{0}, p^{l_2}) = (\det(\mat{R}_2), p^{l_2})$.  Hence
  \begin{align*}
    \eta(\mat{A}_2, \vek{0}, p^{l_3})
    & =
    d_2^2 (\det(\mat{R}_2), p^{l_2}) \\
    & =
    (d_2^2 \det(\mat{R}_2), d_2^2 p^{l_2}) \\
    & =
    (\det(\mat{A}_2), p^{2e_2} p^{l_3 - e_2}) \\
    & =
    (\det(\mat{A}_2), d_2 p^{l_3}).
  \end{align*}
  If $d_2 = p^{l_3}$, then $\mat{A}_2 \equiv \mat{0} \pmod{p^{l_3}}$ and
  $\eta(\mat{A}_2, \vek{0}, p^k) = p^{2l_3}$.  Since $d_2 p^{l_3} = p^{2l_3}$
  it follows that the equality
  \[
    \eta(\mat{A}_2, \vek{0}, p^{l_3}) = (\det(\mat{A}_2), d_2 p^{l_3})
  \]
  also holds in the case when $(\mat{A}_2, p^{l_3}) = p^{l_3}$.  By repeating
  the procedure with Gaussian elimination described before the theorem we get
  a sequence
  \[
    \mat{A}_n, \mat{A}_{n-1}, \ldots, \mat{A}_2
  \]
  of matrices.  Assume that
  \begin{equation}
  \label{eq:inducton-assumption}
    \eta(\mat{A}_m, \vek{0}, p^{l_{m+1}})
    =
    (\det(\mat{A}_m), d_2 d_3^2 \cdots d_m^{m-1} p^{l_{m+1}})
  \end{equation}
  when $m = n - 1, \ldots, 2$.  Let $\mat{R}_n = \mat{A}_n/d_n$, where $d_n =
  p^{e_n} = (\mat{A}_n, p^{l_{n+1}})$.  The exponent of the prime power
  \[
        d_2 d_3^2 \cdots d_{n-1}^{n-2} p^{l_n}
        =
        p^{e_2} p^{2e_3} \cdots p^{(n-2)e_{n-1}} p^{l_n}
  \]
  is
  \begin{multline*}
        (l_3 - l_2) + 2(l_4 - l_3) + \cdots + (n - 2)(l_n - l_{n-1}) + l_n
     \\
        = -l_2 - l_3 - \cdots - l_{n-1} + (n - 1) l_n
        \leq n l_{n+1} = nk,
  \end{multline*}
  since $l_n \leq l_{n+1}$ and all $l_i$ are non-negative integers.
    Hence, by choosing an integer which is an multiplicative inverse to
    $r_{11}$ modulo $p^{nk}$ the equality~\eqref{eq:det-induction-2} will
    be fulfilled when $p^i = d_2 d_3^2 \cdots d_{n-1}^{n-2} p^{l_n}$.  Then
  it follows from~\eqref{eq:eta-induction},~\eqref{eq:det-AnR},~\eqref{eq:det-induction-2}
  and~\eqref{eq:inducton-assumption} that
  \begin{align*}
    \eta(\mat{A}, \vek{0}, p^k)
    & =
    \eta(\mat{A}_n, \vek{0}, p^{l_{n+1}}) \\
    & =
    d_n^n \eta(\mat{A}_{n-1}, \vek{0}, p^{l_n}) \\
    & =
    d_n^n (\det(\mat{A}_{n-1}), d_2 d_3^2 \cdots d_{n-1}^{n-2} p^{l_n}) \\
    & =
    d_n^n (\det(\mat{R}_n), d_2 d_3^2 \cdots d_{n-1}^{n-2} p^{l_n}) \\
    & =
    (d_n^n\det(\mat{R}_n), d_2 d_3^2 \cdots d_{n-1}^{n-2} d_n^n p^{l_n}) \\
    & =
    (\det(\mat{A}_n),
    d_2 d_3^2 \cdots d_{n-1}^{n-2} d_n^{n-1} p^{e_n} p^{l_n}) \\
    & =
    (\det(\mat{A}_n), d_2 d_3^2 \cdots d_n^{n-1} p^{l_{n+1}}) \\
    & =
    (\det(\mat{A}), d_2 d_3^2 \cdots d_n^{n-1} p^k),
  \end{align*}
  since $k = l_{n+1} = l_n + e_n$.

  If $d_i = p^{e_i} = (\mat{A}_i, p^{l_{i+1}}) = p^{l_{i+1}}$ for some $i$,
  then $l_i = 0$ and we will consider the system $\mat{R}_i\vek{x} \equiv
  \vek{0} \pmod{1}$, which have exactly one solution, namely $\vek{x} =
  \vek{0}$.  Further, $d_j = 1$ and $l_j = 0$ for all $j < i$.  The number of
  solutions of the system $\mat{A}_i\vek{x} \equiv \vek{0}
  \pmod{p^{l_{i+1}}}$ is
  \[
    (p^{l_{i+1}})^i
    =
    d_i^i
    =
    d_i^{i-1} p^{e_i}
    =
    d_i^{i-1} p^{l_{i+1}}.
  \]
  Further, we also have that $p^{il_{i+1}} \mid \det(\mat{A}_i)$.  Since $d_2
  d_3^2 \cdots d_{i-1}^{i-2} = 1$ it follows that
  \[
    \eta(\mat{A}_i, \vek{0}, p^{l_{i+1}})
    =
    d_2 d_3^2 \cdots d_i^{i-1} p^{l_{i+1}}.
  \]
  We have proved the theorem.
\end{proof}

\begin{corollary}
  Assume that $p^l$ divide $\det(\mat{A})$ exactly. If $l \leq k$, then
  \[
    \eta(\mat{A}, \vek{0}, p^k) = p^l.
  \]
\end{corollary}

The determinant of a matrix gives an upper bound for how many
solutions there can exist, and therefore the number of solutions will not
increase if we fix the matrix and choose $k$ larger than $l$.  If \mat{A} is not the zero matrix
and $\det(\mat{A}) = 0$, then we have to find all $d_2, d_3, \ldots, d_n$ to determine
the number of solutions.

\begin{example}
  Let
  \[
    \mat{A}
    =
    \begin{pmatrix}
      3 & 6 & 0 \\
      2 & 5 & 1 \\
      6 & 1 & 9
    \end{pmatrix}.
  \]
  Then $\det(\mat{A}) = 60 = 2^2 \cdot 3 \cdot 5$.  Hence, the number of
  solutions of $\mat{A}\vek{x} \equiv \vek{0} \pmod{2^k}$, for $k\geq 2$ is $4$ and for modulo $2$ we have $2$ or $4$ solutions. In this case we get two solutions. For modulo $3^k$ or $5^k$ the number
  of solutions is always $3$ and $5$, respectively.  Note that for all
  other prime numbers $p$ the system has exactly one solution modulo $p^k$.
\end{example}

\begin{example}
  Let
  \[
    \mat{A}
    =
    \begin{pmatrix}
      2 & 2 & 1 \\
      1 & 1 & 2 \\
      1 & 1 & 2
    \end{pmatrix}.
  \]
  Then $\mat{A}\vek{x} \equiv \vek{0} \pmod{3}$ has nine solutions.
\end{example}

\begin{corollary}
  \label{thm:sum-e}
  Let $d_1 = p^{e_1} = (\mat{A}_1, p^{l_2})$.  Then
  \[
    \eta(\mat{A}_n, \vek{0}, p^k)
    =
    d_1 d_2^2 \cdots d_n^n
    =
    p^{\varepsilon},
  \]
  where $\varepsilon = e_1+ 2e_2 + \cdots + n e_n$.
\end{corollary}
\begin{proof}
  If we complete the Gaussian elimination we get the system
  \[
    \begin{ekvsystem}[2]
      r_{11}x_1 &+& r_{12}x_2 &\equiv & 0 \pmod{p^{l_1}} \\
      && a_{22}' x_2 &\equiv & 0 \pmod{p^{l_1}}\rlap{,}
    \end{ekvsystem}
  \]
  where $(r_{11}, p) = 1$ and $\mat{A}_1 = a_{22}'$.  Then $\det(\mat{A}_1) =
  a_{22}'$ and the number of solutions to the system is $d_1 = (a_{22}',
  p^{l_1})$.  From~\eqref{eq:eta-induction} it follows that
  \begin{align*}
    \eta(\mat{A}_n, \vek{0}, p^k)
    & =
    d_2^2 d_3^3 \cdots d_n^n \eta(\mat{A}_1, \vek{0}, p^{l_2}) \\
    & =
    d_1 d_2^2 d_3^3 \cdots d_n^n \\
    & =
    p^{e_1} (p^{e_2})^2 (p^{e_3})^3 \cdots (p^{e_n})^n \\
    & =
    p^{e_1 + 2e_2 + \cdots + n e_n},
  \end{align*}
  and by that the corollary.
\end{proof}

\begin{corollary}
  \label{thm:homogen-m}
  Let $\mat{A}$ be a square integer matrix of order $n$, and $m$ a positive
  integer.  Assume that $m = \prod_{i=1}^N p_i^{k_i}$, where $p_i$ are
  different primes.  Then
  \[
    \eta(\mat{A}, \vek{0}, m)
    =
    \prod_{i=1}^N p_i^{\varepsilon_i}
  \]
  where
  \[
    \varepsilon_i = e_{i,1} + 2e_{i,2} + \cdots + n e_{i,n}
  \]
  is the exponent of the prime power given in Corollary~\ref{thm:sum-e} when
  we consider the system $\mat{A}\vek{x} \equiv \vek{0} \pmod{p^{k_i}}$.
\end{corollary}

\section{Inhomogeneous Systems}
\label{sec:inhomogeneous}

\begin{theorem}
  \label{thm:inhom-system}
  Let $p$ be a prime number, $k$ a positive integer, $\mat{A}_n$ a square integer
  matrix of order $n$, and $\vek{b}_n$ an integer vector of length $n$.  The
  inhomogeous linear system
  \begin{equation}
    \label{eq:inhom-system}
    \mat{A}_n\vek{x} \equiv \vek{b}_n \pmod{p^k}
  \end{equation}
  is solvable if and only if $d_n \mid \vek{b}_n$ and $d_i \mid \vek{b}_i$ for
  all $i = n - 1, \ldots, 2, 1$, where the integers $d_i$ are given in
  similar way as in Theorem~\ref{thm:homogen-system} and the vectors
  $\vek{b}_i$ is the last $i$ elements in the right hand side of the system
  after the $j$th step of the Gaussian elimination.  Moreover, if the system
  is solvable then
  \[
    \eta(\mat{A}_n, \vek{b}_n, p^k)
    =
    \eta(\mat{A}_n, \vek{0}, p^k)
    =
    p^{\varepsilon},
  \]
  where $\varepsilon$ is the exponent given by Corollary~\ref{thm:sum-e}.
\end{theorem} 
\begin{proof}
  The system~\eqref{eq:inhom-system} is solvable when $d_n$ divide each
  component of $\vek{b}$.  Let $\mat{R}_n = \mat{A}_n/d_n$ and $\vek{r}_n =
  \vek{b}_n/d_n$.  Reduce the system in the same way as described
  for~\eqref{eq:system-pl-reduced}.  We get the inhomogeous system
  \[
    \mat{A}_{n-1} \vek{x} \equiv \vek{b}_{n-1} \pmod{p^{l_n}},
  \]
  with $n - 1$ unknowns.  In order for this system to have solutions $d_{n-1}
  = (\mat{A}_{n-1}, p^{l_n})$ must divide all the components of
  $\vek{b}_{n-1}$.  If we continue in this way we get that
  \eqref{eq:inhom-system} is solvable if and only if $d_n \mid \vek{b}$ and
  $d_i \mid \vek{b}_i$ for all $i = n - 1, \ldots, 2, 1$.  The number of
  solutions, if they exists, are the same as in homogeneous case---the
  backward substitution result in the same number of solutions in each step.
\end{proof}

In the following example we will use the algorithm described in Appendix~\ref{sec:algorithm}.

\begin{example}
  Study the linear congruence system
  \begin{equation}
    \label{eq:cs-example}
    \begin{ekvsystem}[5]
      123 x_{1} &+& 152 x_{2} &+& 28 x_{3} &+& 22 x_{4} &+& 144 x_{5} &
      \equiv& 193 \\
      38 x_{1} &+& 189 x_{2} &+& 127 x_{3} &+& 171 x_{4} &+& 141 x_{5} &
      \equiv& 2 \\
      132 x_{1} &+& 232 x_{2} &+& 215 x_{3} &+& 22 x_{4} & &  &
      \equiv& 96 \\
      155 x_{1} &+& 30 x_{2} &+& 178 x_{3} &+& 142 x_{4} &+& 127 x_{5} &
      \equiv& 198 \\
      194 x_{1} &+& 171 x_{2} &+& 16 x_{3} &+& 24 x_{4} &+& 98 x_{5} &
      \equiv& 162
    \end{ekvsystem}
    \pmod{243},
  \end{equation}
  where $243 = 3^5$.  Then $l_6 = 5$.
  Let
  \[
    \mat{A}_5
    =
    \begin{pmatrix}
      123 & 152 & 28 & 22 & 144 \\
      38 & 189 & 127 & 171 & 141 \\
      132 & 232 & 215 & 22 & 0 \\
      155 & 30 & 178 & 142 & 127 \\
      194 & 171 & 16 & 24 & 98
    \end{pmatrix}
    \quad\text{and}\quad
    \vek{b}_5
    =
    \begin{pmatrix}
      193 \\
      2 \\
      96 \\
      198 \\
      162
    \end{pmatrix}.
  \]
  Then $3^2$ is the largest power of $3$ which divide $\det(\mat{A}_5) =
  -134\,741\,218\,779$.  Hence, the number of solutions of the congruence
  system can not exceed~$9$.  We have that $d_5 = (\mat{A}_5, 3^5) = 1$, and
  therefore is $e_5 = 0$, $l_5 = l_6 - e_5 = 5$ and $\mat{R}_5 = \mat{A}_5$.
  Since $d_5$ divides all elements in $\vek{b}_5$ we can continue.  Next, we
  interchange the first and second row.  Hence, $r_{11} = 38$ and
  $r_{11}^{-1} = 32$ modulo $3^5$.  Note that $r_{11}^{-1} =
  200\,673\,618\,026$ modulo $p^{nk} = 3^{5 \cdot 5}$.  After the first step
  in the Gaussian elimination we get the matrices
  \[
    \begin{pmatrix}
      38 & 189 & 127 & 171 & 141 \\
      0 & 71 & 7 & 76 & 180 \\
      0 & 151 & 68 & 157 & 9 \\
      0 & 84 & 114 & 52 & 121 \\
      0 & 63 & 135 & 123 & 56
    \end{pmatrix}
    \quad\text{and}\quad
    \begin{pmatrix}
      2 \\
      97 \\
      153 \\
      241 \\
      139
    \end{pmatrix}.
  \]
  Hence,
  \[
    \mat{A}_4
    =
    \begin{pmatrix}
      71 & 7 & 76 & 180 \\
      151 & 68 & 157 & 9 \\
      84 & 114 & 52 & 121 \\
      63 & 135 & 123 & 56
    \end{pmatrix}
    \quad\text{and}\quad
    \vek{b}_4
    =
    \begin{pmatrix}
      97 \\
      153 \\
      241 \\
      139
    \end{pmatrix}.
  \]
  We have that $d_4 = (\mat{A}_4, 3^{l_5}) = 1$, and therefore is $e_4 = 0$,
  $l_4 = l_5 - e_4 = 5$ and $\mat{R}_5 = \mat{A}_5$.  Since $d_4$ divides all
  elements in $\vek{b}_4$ we can continue.  Further, we do not have to
  interchange any rows or columns this time since $r_{11} = 71$ is relatively
  prime to $3$.  Next step in the Gaussian elimination gives us the matrices
  \[
    \mat{A}_3
    =
    \begin{pmatrix}
      36 & 122 & 54 \\
      27 & 10 & 175 \\
      9 & 213 & 218
    \end{pmatrix}
    \quad\text{and}\quad
    \vek{b}_3
    =
    \begin{pmatrix}
      22 \\
      181 \\
      94
    \end{pmatrix}.
  \]
  We have that $d_3 = (\mat{A}_3, 3^{l_4}) = 1$, and therefore is $e_3 = 0$,
  $l_3 = l_4 - e_3 = 5$ and $\mat{R}_3 = \mat{A}_3$.  Since $d_3$ divides all
  elements in $\vek{b}_3$ we can continue.  Further, we have to interchange
  the first and second column.  Then $r_{11} = 122$.  Next step in the
  Gaussian elimination gives us the matrices
  \[
    \mat{A}_2
    =
    \begin{pmatrix}
      36 & 67 \\
      225 & 56
    \end{pmatrix}
    \quad\text{and}\quad
    \vek{b}_2
    =
    \begin{pmatrix}
      227 \\
      199
    \end{pmatrix}.
  \]
  We have that $d_2 = (\mat{A}_2, 3^{l_3}) = 1$, and therefore is $e_2 = 0$,
  $l_2 = l_3 - e_2 = 5$ and $\mat{R}_2 = \mat{A}_2$.  Since $d_2$ divides all
  elements in $\vek{b}_2$ we can continue.  Further, we have to interchange
  the first and second column.  Then $r_{11} = 67$.  Next step in the
  Gaussian elimination gives us the ``matrices''
  \[
    \mat{A}_1 = 126
    \quad\text{and}\quad
    \vek{b}_1 = 216.
  \]
  If we put these results together we  get the matrices
  \[
    \begin{pmatrix}
      38 & 189 & 171 & 141 & 127 \\
      0 & 71 & 76 & 180 & 7 \\
      0 & 0 & 122 & 54 & 36 \\
      0 & 0 & 0 & 67 & 36 \\
      0 & 0 & 0 & 0 & 126
    \end{pmatrix}
    \quad\text{and}\quad
    \begin{pmatrix}
      2 \\
      97 \\
      22 \\
      227 \\
      216
    \end{pmatrix},
  \]
  which corresponds to the linear congruence system
  \begin{equation}
    \label{eq:trapp}
    \begin{ekvsystem}[5]
      38y_{1} &+& 189y_{2} &+& 171y_{3} &+& 141y_{4} &+& 127y_{5} &
      \equiv& 2 \pmod{3^{5}} \\
      & & 71y_{2} &+& 76y_{3} &+& 180y_{4} &+& 7y_{5} &
      \equiv& 97 \pmod{3^{5}} \\
      & &  & & 122y_{3} &+& 54y_{4} &+& 36y_{5} &
      \equiv& 22 \pmod{3^{5}} \\
      & &  & &  & & 67y_{4} &+& 36y_{5} &
      \equiv& 227 \pmod{3^{5}} \\
      & &  & &  & &  & & 14y_{5} &
      \equiv& 24 \pmod{3^{3}}\rlap{,}
    \end{ekvsystem}
  \end{equation}
  where
  \[
    \begin{ekvsystem}[1]
      38^{-1} &\equiv& 32\pmod{3^{5}} \\
      71^{-1} &\equiv& 89\pmod{3^{5}} \\
      122^{-1} &\equiv& 2\pmod{3^{5}} \\
      67^{-1} &\equiv& 214\pmod{3^{5}} \\
      14^{-1} &\equiv& 2\pmod{3^{3}}
    \end{ekvsystem}
    \quad\text{and}\quad
    \begin{ekvsystem}[1]
      x_{1} &=& y_{\sigma(1)} = y_{1} \\
      x_{2} &=& y_{\sigma(2)} = y_{2} \\
      x_{3} &=& y_{\sigma(3)} = y_{5} \\
      x_{4} &=& y_{\sigma(4)} = y_{3} \\
      x_{5} &=& y_{\sigma(5)} = y_{4}\rlap{.}
    \end{ekvsystem}
  \]
  Where
  \[
    \sigma =
    \begin{pmatrix}
      1 & 2 & 3 & 4 & 5 \\
      1 & 2 & 5 & 3 & 4
    \end{pmatrix}
  \]
  denotes the permutations of the columns done in the computation.
  Since $d_1 = (\mat{A}_1, 3^{l_2}) = (126, 3^5) = 3^2$ we have that $e_1 =
  2$ and $l_1 = l_2 - e_1 = 3$, and therefore is the number of solutions
  \[
    \eta(\mat{A}, \vek{b}, 3^5) = 3^{e_1+2e_2+3e_3+4e_5+5e_5} = 3^{2+0+0+0+0} = 9.
  \]
  Backward substitution gives the following.  First we have that
  \[
    x_3 = y_5 \equiv 2 \cdot 24 \equiv 21 \pmod{3^3}.
  \]
  We lift the result to $\Z_{243}$, that is,
  \[
    x_3 = y_5 = 21 + 3^{l_1} h = 21 + 3^3 h,\qquad h \in \Z_{d_1} = \Z_9.
  \]
  Hence, $x_3 \in \{21, 48, 75, 102, 129, 156, 183, 210, 237\}$.  Since $36$
  is divisible by~$9$, neither $y_4$ or $y_3$ depends on $y_4$, which follows
  from
  \[
    36 \cdot (21 + 3^3 h) \equiv 36 \cdot 21 \equiv 27 \pmod{3^5}.
  \]
  We get that $x_4 = y_3 = 179$ and $x_5 = y_4 = 32$.  The two first
  congruences in the system~\eqref{eq:trapp} gives us the nine solutions:
  \[
    \begin{array}{ccc}
      (43,127,21,179,32) & (151,73,48,179,32) & (16,19,75,179,32) \\
      (124,208,102,179,32) & (2154,129,179,32) & (97,100,156,179,32) \\
      (205,46,183,179,32) & (70,235,210,179,32) & (178,181,237,179,32).
    \end{array}
  \]
  We have solve the system~\eqref{eq:cs-example}.
\end{example}

\section{Discussion}\label{sec:discussion}

The main motivation to investigate the problem was to find a formula. as simple as possible, for the number of solutions of a linear congruence system.

We made the choice to use as elementary methods as
possible. Mainly because the problem is of an elementary nature but also
because of the interest from applications.  It is possible to lift both the
problem and its solution into the context of free modules,
see~\cite{Lang2004}, over the ring $\mathbb{Z}/m\mathbb{Z}$.  The derivation
will when be almost identical but with different vocabulary.

By small changes it is also possible to solve system of congruence equations
with different modulus. We construct a system modulo the least common
multiple of all the moduli, that are equivalent to the original system. This
is the same technique that is mentioned in~\cite{ButSte1955}.

\appendix

\section{Algorithm}
\label{sec:algorithm}

Let $p$ be a prime number, $n$ and $k$ positive integers, $\mat{A} =
(a_{ij})_{n \times n}$ an integer matrix, and $\vek{b} = (b_i)_{n \times 1}$
a vector with integer entries.  The following algorithm determents if the
congruence
\[
  \mat{A}\vek{x} \equiv \vek{b} \pmod{p^k}
\]
is solvable, and in that case it finds the number of solutions, denoted $\eta$,
and the set $X$ of all solutions.
\begin{compactenum}[1.]
\item {[Initialization]} %
  Set $l_{n+1} \gets k$, $s \gets n$ and $t \gets 1$.

\item {[Factor]} %
  \label{step:factor}%
  Set $d_s \gets (\det(\mat{A}), p^{l_{s+1}}) = p^{e_s}$ and $l_s \gets l_{s+1}
  - e_s$

\item {[Solvable?]} %
  If $d_s \centernot{\mid} b_i$ for some $i = t, \ldots, n$, then stop and
  return ``Not solvable''.

\item {[Cancel factor]} %
  Set $a_{ij} \gets a_{ij}/d_s$ and $b_i \gets b_i/d_s$ for all $i, j = t,
  \ldots, n$.

\item {[Zero matrix?]} %
  If $l_s = 0$, then set $e_i \gets 0$ and $l_i \gets 0$ for all $i = s + 1,
  \ldots, n$ and go step~\ref{step:number-of-solutions}.

\item {[Pivoting]} %
  Find an element among $a_{ij}$ for $i, j = t, \ldots, n$ which is relative
  prime to $p$, and perform, if necessary, an interchange of rows or columns
  so that the element is the $t$th element i the main diagonal of \mat{A}.

\item {[Gaussian elimination]} %
  For $i, j = t + 1, \ldots, n$ set $a_{tj} \gets 0$,
  \[
    a_{ij} \gets a_{ij} - a_{tt}^{-1} a_{it} a_{tj} \bmod{p^{nk}}
    \quad\text{and}\quad
    b_i \gets b_i - a_{tt}^{-1} a_{it} b_i \bmod{p^{nk}}.
  \]
  The arithmetic is done modulo~$p^{nk}$ to make sure that the formula in
  step~\ref{step:number-of-solutions} give the correct number of solutions.

\item {[Done?]} %
  If $s > 2$, then set $s \gets s - 1$ and $t \gets t + 1$, and go to
  step~\ref{step:factor}.

\item {[Number of solutions]} %
  \label{step:number-of-solutions}%
  Set $\eta \gets p^{e_1+2e_2+\cdots+ne_n}$.
\end{compactenum}
To determine the solutions the algorithm can be continued in the following
way.
\begin{compactenum}[1.]
\item {[Initialization]} %
  Set $s \gets 2$ and $t \gets n - 1$.

\item {[Introduce solution set]} %
  Let $X$ be the set of all vectors $(\cdot, \ldots, \cdot, y_n)$, where
  $y_n$ is one of the $p^{e_1}$ solutions of $a_{nn} y \equiv b_n
  \pmod{p^{l_1}}$.  If $l_1 \neq l_2$, then for each $(\cdot, \ldots, \cdot,
  y) \in X$ add to $X$ all vectors $(\cdot, \ldots, \cdot, y + i p^{l_1})$,
  for $i = 0, \ldots, p^{e_1} - 1$.

\item {[Backward substitution]} %
  \label{step:BS}%
  For each $\vek{y} = (\cdot, \ldots, \cdot, y_{t+1}, \ldots, y_n) \in X$ set
  \[
    y_t
    \gets
    a_{tt}^{-1}
    \left(
      b_t - \sum_{j=t+1}^n a_{tj} y_j
    \right)
    \bmod{p^{l_s}}
  \]
  and store $y_t$ at position $t$ in \vek{y}.  Note that $(a_{tt}, p) = 1$.

\item {[Lifting]} %
  \label{step:lift}%
  If $l_s \neq l_{s+1}$, then for each $\vek{y} = (\cdot, \ldots, \cdot, y_t,
  \ldots, y_n) \in X$ add to $X$ all vectors $(\cdot, \ldots, \cdot, y_t +
  i_t p^{l_s}, \ldots, y_n + i_n p^{l_s})$, where $i_t, \ldots, i_n = 0,
  \ldots, p^{e_s} - 1$.

\item {[Done?]} %
  If $t > 1$, then set $s \gets s + 1$ and $t \gets t - 1$, and go to
  step~\ref{step:BS}.

\item {[Rearrange solution]} %
  Let $\sigma \in S_n$ be the permutation which summarize the interchanges of
  columns.  Replace each $(y_1, \ldots, y_n) \in X$ with $(x_1, \ldots,
  x_n)$, where $x_i = y_{\sigma(i)}$.
\end{compactenum}

\end{document}